\documentclass[11pt]{article}
\usepackage{amssymb,latexsym,epsfig,epsf,xypic}
\usepackage[mathscr]{euscript}

\makeatletter \@addtoreset{equation}{section} \makeatother

\newtheorem{Lemma}[equation]{Lemma}

\newtheorem{Definition}[equation]{Definition}

\newtheorem{Conj}[equation]{Conjecture}
\newenvironment{Proof}{\noindent\emph{Proof\ }}{\hfill$\square$\\}
\newenvironment{Remarks}{\noindent\textbf{Remarks}\ }

\newcommand{\Rt}[1]{\stackrel{#1\,}{\longrightarrow}}
\newcommand\A{\mathcal A}
\newcommand\db{\bar\partial}
\newcommand\R{\mathbb R}
\newcommand\C{\mathbb C\,}

\newcommand\Pee{\mathbb P}

\newcommand\OO{\mathscr O}

\newcommand\into{\hookrightarrow}
\newcommand\res{\arrowvert_}
\newcommand\tr{\mathrm{tr}}
\newcommand\id{\mathrm{id}}

\newcommand\ip{{\mbox\,\_\hspace{-1.5mm}\shortmid\hspace{1.5pt}}}
\newcommand\ImO{\,\mathrm{Im}\,\Omega}

\title{\textbf{Moment maps, monodromy and mirror manifolds}}
\author{R.\,P. Thomas}
\date{}

\begin{document}

\maketitle
\begin{abstract} \noindent
Via considerations of symplectic reduction, monodromy, mirror
symmetry and Chern-Simons functionals, a conjecture is proposed on
the existence of special Lagrangians in the
hamiltonian deformation class of a given Lagrangian submanifold of a
Calabi-Yau manifold. It involves a stability condition
for graded Lagrangians, and can be proved for the simple case of
$T^2$.
\end{abstract}

\section{Introduction}
Just as explicit solutions of the Einstein and Hermitian-Yang-Mills
equations exist only on spaces that are either low dimensional,
non-compact and/or highly symmetric,
so the equations for special Lagrangian (SLag) cycles,
also important in physics, have the same properties.
Physically there are also similarities in that we have two first order
supersymmetric minimal energy equations (HYMs and SLag)
implying the more standard second order equations (YMs and minimal
volume equations).

There are powerful existence results of Calabi and Yau (and more
recently Tian, Donaldson and others) for the Einstein equations, and of
Donaldson-Uhlenbeck-Yau for the HYM equations, so long as we are on a
K\"ahler (or projective) manifold; this often reduces
an infinite dimensional problem in PDEs to a finite dimensional
problem in linear algebra. Producing many K\"ahler-Einstein
(e.g. Calabi-Yau) manifolds becomes trivial, and dealing with
Hermitian-Yang-Mills connections requires only algebraic
computations; in both cases the complicated role of the K\"ahler form
and/or metric is almost removed. This can be thought of as possible
because of the existence of some infinite dimensional geometry recasting
the equations in terms of moment maps and symplectic reduction. A
similar situation for SLags would therefore be highly desirable.
In particular it might give a way of studying SLags using only
Lagrangians and symplectic geometry, much as HYM connections are
studied via stable bundles and algebraic geometry.

This paper explores the mirror symmetry of holomorphic bundles
(on a Calabi-Yau $3$-fold $M$,
often referred to here as `the complex side') and Lagrangians
(on the mirror Calabi-Yau 3-fold $W$, 'the symplectic side', known
as the K\"ahler side in the physics literature). Many
people have worked and are still working on proving some kind of
direct correspondence between such objects given an
SYZ torus fibration \cite{SYZ}; see for example \cite{AP}, \cite{BBHM},
\cite{Ch}, \cite{Fu1}, \cite{Gr}, \cite{LYZ}, \cite{PZ}, \cite{Ty}, and
see \cite{MMM} for a review of this and many many more issues in mirror
symmetry. Here, however, we work purely formally without
reference to a particular pair of mirror manifolds, without
worrying about what mirror symmetry might rigorously mean, and we will
not try to give any explicit correspondence. Using mirror symmetry
merely as motivation, we point out some similar structures on
both sides of the mirror map. Under some conditions (in some
`large complex structure' or `semi-classical' or somesuch limit)
these structures might be genuinely dual; again it does not matter
if they are not in general. For instance, physics \cite{MMMS},
\cite{DFR} predicts that one should consider not the HYM equations and
slope but some perturbation of them away from the large complex
structure limit; however these equations also come from a moment map
and, conjecturally, a stability condition (for a discussion of such
matters see \cite{Le} or \cite{T3}). So while the slope and phase of
Lagrangians discussed below might not be exactly mirror to slope of
bundles, it should be mirror to something with analogous properties
and significance.

Loosely, we would like to think of submanifolds in a fixed homology class
as mirror to connections on a fixed
topological complex bundle (with Chern classes mirror to the homology
class); then Lagrangians should correspond to holomorphic
connections (i.e. integrable connections; those with no
$(0,2)$-curvature) and special Lagrangians to those with HYM
curvature. These last two conditions should be stability conditions
for the group actions of hamiltonian deformations and complex gauge
transformations, respectively. The full picture is much more complicated,
involving triangulated categories and so forth, as envisaged some six
years ago in the seminal conjecture of Kontsevich \cite{K}; we can
ignore this in only using mirror symmetry
as motivation. It could be noted, however, that the functionals defined
below are additive under exact sequences of holomorphic vector bundles
and sums of Lagrangians, so should extend to the derived category
of coherent sheaves and the derived Fukaya category of Lagrangians
respectively.

First note that while the connections side has a complex
structure and a complex gauge group involved, the Lagrangian side needs
complexifying. So motivated by Kontsevich \cite{K} and by physics
(e.g. \cite{SYZ}) we add in connections on the submanifolds (which
will later reduce to flat connections on Lagrangians). The dictionary
we are aiming towards, much of which is already standard, is the
following in the 3-dimensional case; all the terms used will be defined
in due course. \newpage

\renewcommand\arraystretch{2}
$$
\begin{array}{|c|c|}
\hline \text{Complex side }M & \text{Symplectic side }W \\
\hline\hline \Omega=\Omega^{ }_M\in H^{3,0} & \omega=\omega^{ }_W\in
H^{1,1} \\ \hline
H^{\mathrm{ev}} & H^3 \\ \hline
\mathrm{Connections\ }A\mathrm{\ on\ a\ fixed\ }C^\infty\mathrm{\ complex}
& \mathrm{Submanifolds/cycles\ }L\mathrm{\ in\ a\ fixed
\ class} \vspace{-4mm} \\
\mathrm{bundle\ }E;\ \,v(E)=ch(E)\sqrt{\mathrm{Td\,}X}\in H^{\mathrm{ev}} &
[L]\in H^3,\mathrm{\ with\ a\ connection \ on\ }\C\times L \\ \hline
CS_\C(A=A_0+a)= \hspace{25mm} &
f_\C(A,L)=\int_{L_0}^L(F+\omega)^2 \hspace{2cm}\vspace{-3mm} \\
\hspace{6mm} \frac1{4\pi^2}\int_M\tr\left(\db_{A_0}a\wedge
a+\frac23a^{\wedge3}\right)\wedge\Omega &
\hspace{1cm}=\int_{L_0}^L(F^2+\omega^2)+2
\int_{L_0}^L\omega\wedge F \\ \hline
\mathrm{Critical\ points:\ }F_A^{0,2}=0 &
\mathrm{Critical\ points:\ }\omega\res L=0,\ F_A=0 \vspace{-4mm} \\
\mathrm{holomorphic\ bundles\ } &
\mathrm {Lagrangians\ +\ flat\ line\ bundles\ } \\ \hline
\text{Holomorphic Casson invariant \cite{T1}} &
\text{Counting SLags \cite{J}} \\ \hline
\mathrm{Gauge\ group} & U(1)\mathrm{\ gauge\ group\ on\ }L \\ \hline
\mathrm{Complexified\ gauge\ group} & \mathrm{Hamiltonian\ 
deformations} \\ \hline
\omega=\omega^{\ }_M\in H^{1,1} & \Omega=\Omega^{\ }_W\in H^{3,0} \\ \hline
\mathrm{Moment\ map\ \ }F_A\wedge\omega^{n-1} &
\mathrm{Moment\ map\ \,}\ImO\res L \\ \hline
\mathrm{Stability,\ slope\ }\mu=\frac1{\mathrm{rk\,}E}
\int\tr\,F_A\wedge\omega^{n-1} & \mathrm{Stability,\ slope\ }
\mu=\frac1{\mathrm{vol\,}(L)}\int_L\ImO \\ \hline
\end{array}
$$
(In the fourth line, $v(E)$ is the Mukai vector of $E$; in the last
line, vol$\,(L)$ is the cohomological volume measured with respect
to Re$\,\Omega$.) A SLag cycle (of phase $\phi$) in a Calabi-Yau
is a Lagrangian with Im\,$(e^{-i\phi}\Omega)\res L\equiv0$ \cite{HL};
then Re\,$(e^{-i\phi}\Omega)\res
L$ is the Riemannian volume form on $L$ induced by the Ricci-flat
metric on the Calabi-Yau. Obviously, rotating $\Omega$ by $e^{-i\phi}$
gives SLags in the more traditional sense of phase zero.
The part of the theory to do with SLags will apply in all dimensions;
it is only the functionals that are special to Calabi-Yau 3-folds.

We will partially justify the above table, though the symplectic
structure and moment map give problems that will appear in due course.
However we can derive enough to arrive at
a conjecture about Lagrangians and SLags for which evidence will
be given by using monodromy and mirror symmetry from \cite{ST}
to interpret an example of Lawlor and Joyce \cite{J}. \\

\noindent \textbf{Acknowledgements.} The debt of any ex-student of
Simon Donaldson writing a paper on moment maps should be clear.
This work is also more immediately influenced by the papers
\cite{D},\,\cite{J},\,\cite{K}. In particular I was surprised to see the
Lagrangian condition coming from a moment map in \cite{D},\,\cite{H},
which does not fit into the
scheme I always supposed was true. So the purpose of this paper, apart
from trying to set a record for the number of m's in a title, is to
expand on that scheme and to try to get the \emph{special} condition
from a moment map instead. This paper was finished in the summer of 2000
and reported on in \cite{T3}; since then exciting new ideas have
appeared in physics \cite{Do} and mathematics \cite{KS} better explaining
mirror symmetry. I would like to thank S.-T. Yau, C. H.
Taubes and Harvard University for support, and Yi Hu,
Albrecht Klemm, Ivan Smith and Xiao Wei Wang for
useful conversations. Communications with Mike Douglas, Dominic Joyce,
Paul Seidel, S.-T. Yau and Eric Zaslow have been extremely influential.

\section{Chern-Simons-type functionals and critical points}

Consider the space $\A$ of $(0,1)$-connections $A$ on a fixed complex
bundle on a Calabi-Yau 3-fold $M$. This infinite dimensional space
has a natural complex structure, with respect to which it admits a
holomorphic functional, Witten's holomorphic Chern-Simons functional
\cite{W1},\,\cite{DT},
$$
CS_\C(A=A_0+a)=\frac1{4\pi^2}\int_M\tr\left(\db_{A_0}a\wedge
a+\frac23a\wedge a\wedge a\right)\wedge\Omega,
$$
where $\Omega$ is the holomorphic (3,0)-form. It is infinitesimally
gauge-invariant (gauge transformations not homotopic to the identity
can give periods to $CS_\C$) and its gradient is $F_A^{0,2}$,
with zeros the integrable connections. That is, after dividing by
gauge equivalence (under which grad$\,CS_\C$ is invariant), the critical
points of $CS_\C$ form the space of holomorphic bundles of the same
topological type. As critical points of a functional, moduli of
holomorphic bundles have virtual dimension zero, and one might try to
make sense of counting them -- a holomorphic Casson invariant \cite{T1}.
This is independent of deformations of the complex structure, but can
have wall-crossing changes as the K\"ahler form varies. (This is because
we count only stable bundles, and the notion of stability depends on a
K\"ahler form.)

On the other hand, on a different Calabi-Yau 3-fold $W$ (for instance
the mirror, in some situation where this makes sense),
Lagrangians are the critical points of a functional
too, on the space of all 3-dimensional submanifolds (or cycles):
$$
f_\R(L)=\int_{L_0}^L\omega\wedge\omega,
$$
where $\omega$ is the symplectic form on $W$. Here $L_0$ is a fixed cycle
in the same homology class, and we integrate over a 4-chain
with boundary $L-L_0$; the functional $f_\R$ is invariant under the choice
of different, homologous, 4-chains (picking non-homologous 4-chains can give
periods to $f_\R$). It is invariant under
deformations of $L$ pulled back from hamiltonian deformations of $W$
(deformations generated by vector fields $v$ on $W$ whose
contraction with $\omega$ is exact $v\ip\omega=dh$ at each point in
time) as $\int_L\omega\wedge dh=0$, and its
gradient is $\omega|_L$. Thus its critical points are Lagrangian
submanifolds. We would like to think of $f_\R$ as mirror
to $CS_\C$, but to do so we must complexify it.

Thus we work on the space $\A$ of submanifolds $L$ of $W$ with $U(1)$
connections
$A$ on the trivial bundle $\C\times L$ on $L$.  Notice these submanifolds
are not parameterised by a map of a real 3-manifold into $W$; we are only
interested in the image $L$. \emph{From now on we shall restrict
attention to smooth Lagrangian submanifolds.} Formally, we
consider the tangent space to $\A$ at a point $(A,L\subset W)$ to be
\begin{equation} \label{tsp}
\Omega^1(L;\R)\oplus \Omega^1(L;\R),
\end{equation}
at least for those $L$ with no $J$-invariant subspaces of its tangent
spaces ($J$ is the complex structure on $W$, and this is reasonable since
we are looking for Lagrangian submanifolds after all). The first factor
is the obvious tangent space to the connections on $L$; the second gives
deformations of $L$ via the vector fields produced by contracting with
the K\"ahler form $\omega$ on $W$. That is, we use the metric on $W$ to
map $\Omega^1(L)$ to $\Omega^1(W)|_L$, then use the isomorphism
provided by $\omega$ to get a vector field along $L$. Equivalently,
using the metric on $W$, we may think of one-forms on $L$ as tangent
vectors to $L$, then apply the complex structure $J$ on $W$ to give
$W$-vector fields on $L$. We denote this map from one-forms to normal
vector fields by
\begin{equation}\label{ip}
\Omega^1(L)\to TW\res L, \qquad \sigma\mapsto\sigma\ip\omega^{-1}.
\end{equation}
Connections on $L$ are carried along by the
vector field to connections on nearby cycles, and we are identifying
the space of $U(1)$ connections with $i\Omega^1(L;\R)$.

There is a natural almost complex structure on $\A$, acting as
\renewcommand\arraystretch{1}
$$
J=\left(\begin{array}{cc} 0 & 1 \\ \!-1 & 0 \end{array}\right),
$$
with respect to the splitting (\ref{tsp}) of the tangent spaces. With
respect to this we claim to have the following holomorphic functional
$$
f_\C(A,L)=\int_{L_0}^L(F+\omega)^2=\int_{L_0}^L(F^2+\omega^2)\ +\ 2
\int_{L_0}^L\omega\wedge F.
$$
Here we have extended $A$ to a connection on the trivial bundle on the
whole of $W$ (restricting to a fixed connection $A_0$ on $L_0$, and to
$A$ on $L$) and taken its curvature form $F$. We have again picked
a 4-cycle bounding $L-L_0$; because $F$ and $\omega$ are closed the
resulting functional is independent of different homologous choices
of the 4-cycle, and in general well defined up to the addition of some
discrete periods. It is also (again) independent of hamiltonian isotopies
of $L$. Notice that the $\int_{L_0}^LF^2$ term is just the
\emph{real} Chern-Simons functional $CS_\R$ of the connection $A$ on $L$,
whose critical points are well known to be flat connections. As pointed
out to me by Eric Zaslow, the real and complex Chern-Simons
functionals already appear in \cite{W1} and \cite{Va} as possible
mirror partners (this is partially justified in \cite{LYZ}), but without
the terms in the symplectic form (and including instanton corrections
from holomorphic discs which we are ignoring for our rough analogy).
Asking for a real function to be equal to a complex one is possible
when one restricts attention to a real slice such as the space of
\emph{Lagrangian} submanifolds in $\A$; deforming within this space
the imaginary part of $f_\C$ remains constant and it reduces to
$CS_\R$. But allowing the imaginary counterparts to these real
deformations the right functional to consider is $f_\C$. Notice also
that if $\omega/2\pi$ is integral, so that we can pick a connection $B$
with curvature $-i\omega$, then the action functional
can be written in the more familiar looking Chern-Simons form
$$
f_\C(A,L)=\int_L(B+iA)\wedge d(B+iA)=\int_LCdC
$$
for the `complexified connection' $C=B+iA$ (a $\C^{\!\times}$-connection,
instead of a $U(1)$-connection.) This makes more contact with the
physics literature and allows one to extend the identification of
$CS_\R$ and $CS_\C$ in \cite{LYZ} to non Lagrangian sections, giving
complex values. Tian has informed me that he and Chen have also
considered the functional $f_\C(A,L)$ \cite{Ch}.

Mirror symmetry should relate Lagrangians not just to bundles but the
whole derived category. For Riemann surfaces $C\subset M$, for instance,
there is a functional in \cite{DT}, \cite{W2} rather like $f_\R$
above:
$$
\int_{C_0}^C\Omega
$$
is formally holomorphic and has as critical points the \emph{holomorphic}
curves $C$. Similarly for four-manifolds $S\subset M$ with connections on
them the following functional (formally similar to $f_\C$)
$$
\int_{S_0}^S\tr\,F\wedge\Omega
$$
has critical points the holomorphic surfaces with flat connection on them.
Alternatively, as $CS_\C$ is additive under extensions of
bundles it does extend to the derived category. (Whether these two
approaches are compatible; i.e. whether or not the functional
associated to a curve or surface is the same as $CS_\C$ applied
to a locally free resolution of its structure sheaf, up to a constant,
seems to not have been worked out.)

That $f_\C$ is holomorphic follows from the computation that the derivative
of $f_\C$ down $a\in\Omega^1(L)\oplus0$ (that only changes the connection
$A\mapsto A+\delta a$) is $\int_L2F\wedge ia+2\omega\wedge ia$,
while the derivative down $-Ja\in0\oplus\Omega^1(L)$, i.e. down the
vector field $a\ip\omega^{-1}$, is $\int_L2\omega\wedge a+2a\wedge F$.
The second expression is $-i$ times
the first, so the derivative is complex linear and $f$ is holomorphic.
Equivalently we are saying that $df_\C$ is the 2-form
$$
2i\,(F+\omega)\ \oplus\ 2\,(F+\omega),
$$
which pairs with the tangent space (\ref{tsp}) by integration over $L$ to
give a form of type (1,0) on (\ref{tsp}).

Thus critical points of the functional are Lagrangian cycles with flat line
bundles on them: exactly the basic building blocks of the objects proposed
in \cite{K} to be mirror dual to the holomorphic bundles that are the
critical points of $CS$. So this ties
in three well known moduli problems of virtual dimension zero (i.e. with
deformation theories whose Euler characteristic vanishes) -- flat bundles
on 3-manifolds, holomorphic bundles on Calabi-Yau 3-folds, and Lagrangians
(up to hamiltonian deformation) in symplectic 6-manifolds.

So as mirror to \cite{T1} one would like to count Lagrangians
(up to hamiltonian deformations) plus flat line bundles on them, and
this is what Joyce's work \cite{J} has begun to tackle (in the rigid
case of $L$ being a homology sphere). Mirroring precisely the behaviour
of the holomorphic Casson invariant this count appears to be independent
of deformations of the K\"ahler form and to have wall-crossing changes
as the complex structure varies.

\section{Gauge equivalence and moment maps}

In fact what Joyce is proposing to count is \emph{special} Lagrangian
spheres with flat line bundles on them
(hence the otherwise anomalous dependence on the complex
structure), while \cite{T1} counts \emph{stable} bundles (i.e. by
Donaldson-Uhlenbeck-Yau, modulo the technicalities of polystable and
non-locally-free sheaves, we count \emph{Hermitian-Yang-Mills}
connections; hence the dependence on the K\"ahler form). (Tyurin
\cite{Ty} was perhaps the first to suggest
that the holomorphic Casson invariant should be related by mirror
symmetry to the real Casson invariant (here the $U(1)$ Casson
invariant) of SLag submanifolds.)

The link should be, of course, that we want to consider holomorphic
connections on one side, up to complex gauge equivalence, and
Lagrangians on the other side, up to hamiltonian isotopy, and in both cases
we try to do this by picking distinguished representatives of equivalence
classes by the usual method of symplectic reduction.
Bringing in a K\"ahler structure on the complex side, we get a moment
map for the gauge group action, whose zeros give the HYM equations.
Dually, we would like to bring in the holomorphic 3-form on the
symplectic (K\"ahler) side, and get a complex group to act.
So again complexify by adding flat line bundles: consider the critical
points of the functional $f$ of the last section, i.e.
the space
$$
\mathcal Z=\{(L,A)\,:\,L\subset W\mathrm{\ is\ Lagrangian,\ }A
\mathrm{\ is\ a\ flat\ connection\ on\ }L\,\}
$$
(\emph{not} up to gauge equivalence). In fact consider this space
on a Calabi-Yau manifold $W$ of any dimension $n$. It has tangent space
$$
T_{(L,A)}\mathcal Z=Z^1(L)\oplus Z^1(L)
$$
($Z^1(L)$ denotes closed real one-forms on $L$), the first being tangent
to the space of flat connections, the second giving
normal vector fields (by contracting with $\omega^{-1}$) preserving the
Lagrangian condition. We have an obvious almost complex structure
\begin{equation} \label{aha}
J=\left(\begin{array}{cc} 0 & 1 \\ -1 & 0 \end{array}\right).
\end{equation}

Then the real group $C^\infty(L;\R)/\R$ acts as the Lie algebra to the
group of gauge transformations on
the flat line bundles (taking $d$ and adding to the connection)
whose complexification $C^\infty(L;\C)/\C$
acts complex linearly: the imaginary part $C^\infty(L;\R)/\R$
acts by hamiltonian deformations through the normal vector field
$$
h\mapsto dh\ip\omega^{-1}.
$$
Unfortunately, without using a metric this vector field is only
defined up to the addition of tangent vector fields to $L$; the map
(\ref{ip}) is really a map to $(TW\res L)/TL$ which we have lifted to
$TW\res L$ using the metric. (Equivalently we can extend
$h$ to a first formal neighbourhood of $L$ in different ways to get a
different vector field.) How we pick this alters how we carry the flat
connection along with $L$, and how the almost complex structure
(\ref{aha}) acts. For instance suppose we are in the rather artificial
case of $L$ being transverse to an SYZ $T^n$-fibration. Then
we can carry $L$ and the flat connection up the fibres
and identify the functions $C^\infty(L)$ from Lagrangian
to Lagrangian using the projection. Thus the group remains constant
as $L$ moves (effectively what we are doing is extending
functions from $L$ to a neighbourhood of $L$ in $W$ by pulling up
along the SYZ fibres). This does not work when $L$ branches over the
base of such a fibration. One can instead use the metric to define
normal vector fields, but then identifying the Lie algebra
$C^\infty(L)$ with a fixed $C^\infty(L_0)$ for all $L$ becomes difficult.

This problem is perhaps not so surprising -- the moment the Lagrangian
has branching over the base of an SYZ fibration simple explicit
correspondences between Lagrangians and vector bundles (such as \cite{LYZ})
also break down due to our ignoring important holomorphic
disc instanton corrections that appear in the physics. For instance
recent work of Fukaya, Oh, Ohta and Ono \cite{FO3}, surveyed in \cite{Fu1},
show these provide the obstructions mirror to those of deformations
of holomorphic bundles \cite{T2} -- one should not in general consider
all (S)Lags (which are unobstructed) as mirror to holomorphic bundles,
but only those whose Floer cohomology (whose definition involves holomorphic
discs) is well defined.

However, what is clear is the totality of the group action,
even if identifying individual elements causes problems, and this is
all we really need. For instance in the $K3$ (or $T^4$) case one can get
the same total group orbit, with a genuine
fixed group acting, by hyperk\"ahler rotating a construction due to
Donaldson \cite{D}. The end result is that
one considers parametrised Lagrangian embeddings $f$ from a Riemann
surface $L$ into the $K3$ such that the pullback of Re$\,\Omega$ is a
fixed symplectic form on $L$. Then the group of exact
symplectomorphisms of
$$
(L,f^*\mathrm{Re\,}\Omega)
$$
provide a symmetry group of the space of maps $f$, which also carries
a natural K\"ahler structure. Complexified orbits give hamiltonian
deformations, and the moment map is $m(f)=f^*\ImO$. The connection
with our construction is that after fixing a line bundle $\eta$
and connection with curvature $f^*\mathrm{Re\,}\Omega$, an infinitesimal
symplectomorphism $\phi$ induces a flat connection, via parallel transport
and pull back, on the bundle $\eta\otimes\phi^*\eta^*$. Globally the
action is different (this action has non-zero Lie bracket, for instance,
and a fixed group) but the total group orbit and the moment map
(see below) are the same.

In general it is clear that the problem of identifying the group
for different embeddings of $L$ should be resolved by working with
the space of maps from a fixed $L_0$ to $W$, and enlarging the
group by including diffeomorphisms of $L_0$, giving a semi-direct
product of Diff$(L_0)$ and $U(1)$ gauge transformations on $L_0$.
Then the moment map for the diffeomorphism
part of the total group would be the Lagrangian condition as in
\cite{D}, and the problems we are encountering would come from
the fact that the group is a semi-direct product and not a product,
so that we cannot separate the two out and divide by them separately,
as in effect we have been trying to do. Unfortunately, I have not found the
correct formulation of the problem, but it is not so important for
follows.

So we shall not worry too much about whether the complex structure
defined above is integrable, the group is fixed, or the symplectic
structure below is closed. In 1 complex dimension it is trivial, in
2 we can use Donaldson's picture, and in 3 dimensions
we could either try to use an abstract SYZ fibration to deform and
identify Lagrangians transverse to it, or take everything
in this section as motivation for finding the stability condition for
Lagrangians of the next section.

Fix a homology class of Lagrangians and multiply $\Omega$ by a unit
norm complex number so that $\int_L\ImO=0$.
We induce a symplectic structure on $\mathcal Z$ from $J$ and the
following metric on the tangent space
$$
\langle a,b\rangle=\int_La\wedge((b\ip\omega^{-1})\ip\ImO),
$$
for $a,\,b$ closed 1-forms. A computation in local coordinates
shows this is symmetric in $a$ and $b$; in fact it can be written as
\begin{equation} \label{lc}
\int_La\wedge(\widetilde b\ip\mathrm{\,Re\,}\Omega)=\int_L\cos\theta
\,(a\wedge*b),
\end{equation}
where $\,\widetilde{\ }\,$ is the isomorphism $T^*L\to TL$ set up by
the induced metric on $L$, $\Omega\res L=e^{i\theta}\mathrm{vol}_L$,
and vol$_L$ the Riemannian volume form on $L$ induced by the Ricci-flat
metric. Thus for Lagrangians with $\theta\in(-\pi/2,\pi/2)$, i.e. 
those for which Re$\,\Omega$ restricts to a nowhere vanishing volume
form on $L$ and so are not too far from being SLag ($\theta\equiv0$),
this gives a non-degenerate metric.

The symplectic form is invariant under the group action, and formally
the moment map is indeed $m(L,A)=\ImO$ in the dual $\Omega^n(L)_0$ of
the Lie algebra (i.e. $n$-forms on $L$ with integral zero). This follows
from the computation
$$
X\int_Lh\ImO=\int_Lh\,d\,(X\ip\ImO)=
\int_Ldh\wedge((\sigma\ip\omega^{-1})\ip\ImO),
$$
where $X=\sigma\ip\omega^{-1}$ is a normal vector field to the Lagrangian
$L$ down which we compute the derivative of the hamiltonian $\int_L
h\ImO=\langle m(L,A),h\rangle$ for the infinitesimal
action of $h$. Here have extended $h$ to a first-order neighbourhood
of $L\subset W$ so that it is constant in the direction of
$X=\sigma\ip\omega^{-1}$. Then the right hand side of the above
equation is the pairing using the symplectic form of $dh$ and $\sigma$,
as required.

Infinitesimally we can see the moment map interpretation very
easily, and fitting naturally with the mirror bundle point of view.
Deformations of holomorphic connections $A$ modulo
complex gauge equivalence are given by a ker$\,\db_A$/im$\,\db_A$\,
first cohomology group, related to deformations
ker$\,\db_A\,\cap$\,ker$\,\db_A^*$ of the HYM equations
(modulo \emph{unitary} gauge transformations)
via Hodge theory, with the moment
map equation providing the $d^*=0$ slice to the imaginary part of the
linearised group action. Similarly, deformations of
Lagrangians are given by closed 1-forms
ker$\,d:\,\Omega^1(L;\R)\to\Omega^2(L;\R)$, so that dividing by hamiltonian
deformations we get
$$
H^1(L)=\mathrm{ker}\,d/\mathrm{im}\,d.
$$
If instead of dividing we impose the \emph{special}
condition, we get a ker$\,d^*$ slice
$$
H^1(L)=\mathrm{ker}\,d\cap\mathrm{ker}\,d^*,
$$
to the (imaginary) deformations (real deformations are given by
changing the flat $U(1)$ connection that can be incorporated into this).

\subsection*{A symplectic example}

To motivate a guess at the correct definition of stability for
Lagrangians, we expand on an example of Lawlor and Joyce
(\cite{J} Sections 6 and 7, building on work of \cite{Ha}, \cite{L};
see also a similar example in \cite{SV} that is studied in
\cite{TY}), explaining its relevance to mirror symmetry, and giving a
simple example in algebraic geometry that mirrors it.

First define the pointwise phase $\theta$ of a submanifold $L$: we may
write
$$
\Omega\res L=e^{i\theta}\mathrm{\,vol\,}
$$
where vol is the Riemannian volume form on $L$ induced by Yau's
Ricci-flat metric \cite{Y} on $W$. Thus vol provides a (local) orientation
for $L$, and reversing its sign alters the phase $\theta$ by $\pi$. A
SLag is a Lagrangian with \emph{constant} phase $\theta$.

At first sight $\theta$ is multiply-valued; we always choose it to be a 
fixed single-valued function to $\R$, lifting $e^{i\theta}:\,L\to S^1$
and thus providing the Lagrangian with a \emph{grading} as introduced by
Kontsevich \cite{K},\,\cite{S2}. Thus \emph{we only consider Lagrangians
of vanishing Maslov class} -- for a Calabi-Yau this is the winding class
$\pi_1(L)\to\pi_1(S^1)$ of the phase map
$$
L\Rt{e^{i\theta}}S^1,
$$
which of course vanishes for a SLag. (The definition of grading in
\cite{K}, \cite{S2} is topological and uses the universal
$\mathbb Z$-cover of the bundle of Lagrangian Grassmannians; here we
first pass to the $\mathbb Z/2$ orientation cover of the Grassmannian,
choosing an orientation of our Lagrangians,
and then use a complex structure to pass to the universal
$\mathbb Z$-cover of this. The two definitions
are of course equivalent.)

Similarly we can define a kind of average phase $\phi=\phi(L)$ of a
submanifold (or homology class) $L\subset W$ by
$$
\int_L\Omega=A\,e^{i\phi(L)},
$$
for some real number $A$; we then use Re\,$(e^{-i\phi(L)}\Omega\res L)$
to orient $L$. Reversing the sign of $A$ alters the phase by $\pi$
and reverses the orientation. Again for a \emph{graded} Lagrangian
$L=(L,\theta)$, and \emph{we will always implicitly assume a grading},
$\phi(L)$ is canonically a real number (rather than $S^1$-valued).
Shifting the grading $[\,2n\,]:\,\theta\mapsto\theta+2n\pi$ gives a
similar shift to the phase $\phi(L)$.

The terminology comes from the fact that if there is a submanifold
in the same homology class as $L$ that is SLag with respect to some
rotation of $\Omega$, then it is with respect to $e^{-i\phi(L)}\Omega$.
Slope, which we define as
$$
\mu(L):=\tan(\phi(L))=\frac1{\int_L\mathrm{Re\,}\Omega}\int_L\ImO,
$$
is defined independently of grading, is monotonic in $\phi$ in the
range $(-\pi/2,\pi/2)$, and is invariant under change of orientation
$\phi\mapsto\phi\pm\pi$. This agrees with the slope of a straight
line SLag in the case of $T^2$, as featured in \cite{PZ}, and we think
of it as mirror to the slope of a mirror sheaf, as is shown for tori
in \cite{PZ} (see \cite{DFR} for corrections in higher dimensions away
from the large complex structure limit).

Joyce describes examples of SLags which we interpret as follows.
We have a family of Calabi-Yau 3-folds $W^t$ as $t$ ranges through (a
small open subset of) the moduli space of complex structures on $W$
\emph{with fixed symplectic structure}. That is, the holomorphic 3-form
$\Omega^t$ varies with $t$, but the K\"ahler form $\omega$ is fixed.
We also have a family of SLag homology 3-spheres
$L_1^t,\ L_2^t\subset W^t$ such that $L^t_1$ and $L^2_t$ intersect
at a point. If we choose a rotation of $\Omega^t$ such that $L_2^t$
always has phase $\phi_2^t\equiv0$ (this is possible locally at least;
in the family described later it will have to be modified slightly),
then we are interested as $t$ varies only in the complex number
$$
\int_{L_1^t}\Omega^t=R^te^{\phi^t_1}
$$
and its polar phase $\phi=\phi^t_1$; we plot this (i.e. the projection
from the complex structure moduli space to $\C$ via this map) in Figure 1.

Then in Joyce's example, for $\phi<0$ (and $R^t>0$)
there is a SLag $L^t$ (of some phase $\phi^t$)
in the homology class $[L^t]=[L^t_1]+[L^t_2]$,
such that as $\phi\uparrow0$, this degenerates to a singular union
of SLags of the same phase $L^t=L^t_1\cup L^t_2$ and then disappears
for $\phi>0$.

Most importantly, where $L^t$ exists as a smooth SLag ($\phi<0$) we have
the slope (and phase) inequality
\begin{equation} \label{ineq}
\mu^t_1<\mu^t_2, \qquad \mathrm{i.e.\ \,}\phi=\phi^t_1<\phi^t_2\equiv0;
\end{equation}
at $t=0,\ L^t$ becomes the singular union of $L^t_1$ and $L^t_2$,
with
$$
\mu^t_1=\mu^t_2\ (\phi^t_1=\phi^t_2);
$$
then there is no SLag in $L$'s homology class for
$$
\mu^t_1>\mu^t_2\ (\phi^t_1>\phi^t_2),
$$
\emph{though there is a Lagrangian}, of course -- the symplectic
structure has not changed. Though
we have been using slope $\mu$ in order to strengthen the analogy
with the mirror (bundle) situation, from now on we shall use only
the phase (lifted to $\R$ using the grading). While each is
monotonic in the other for small phase (as $\tan\phi=\mu$), slope
does not see orientation as phase does; reversing orientation adds
$\pm\pi$ to the phase but leaves $\mu$ unchanged. This is related to
the fact that we should really be working with complexes and so forth
on the mirror side (the bundle analogy is too narrow) and changing
orientation has no mirror analogue in terms of only stable bundles; it
corresponds to shifting (complexes of) bundles by one place
in the derived category. While slopes of bundles cannot go past
infinity (without moving degree in the derived category at least), for
Lagrangians they certainly can, and phase $\phi$ continues monotonically
upwards as its slope $\tan\phi$ becomes singular and then negative.

Importantly, we can think of the various SLags as independent of time
when thought of as Lagrangians in the fixed symplectic manifold $W^t$:

\begin{Lemma}
For $t>0$ the SLags $L^t$ are all in the same hamiltonian deformation
class. Similarly for $L_1^t,\ L^t_2$, and for $t<0$.
\end{Lemma}

\begin{Proof}
Now choosing the phase of $\Omega^t$ such that $\phi(L^t)\equiv0$,
\begin{equation} \label{fred}
\int_L\frac{d}{dt}(\ImO^t)=\int_L\,\mathrm{Im}\,\dot\Omega^t=0.
\end{equation}
To show this deformation preserves the hamiltonian class of L,
we need to find a corresponding first order hamiltonian deformation
$dh\ip\omega^{-1}$ under which the change in $\ImO^t$,
$$
\mathcal L_{dh\_\hspace{-.4mm}\shortmid\hspace{1pt}\omega^{-1}}
(\ImO^t)\res L=d((dh\ip\omega^{-1})\ip\ImO^t)\res L,
$$
is $-$Im\,$\dot\Omega^t\res L$. But as Re\,$\Omega^t\res L$ is the induced
Riemannian volume form vol$^t$ on $L$, this means we want to solve
$$
-\mathrm{Im\,}\dot\Omega^t\res L=d(J(dh\ip\omega)\ip\mathrm{\,Re\,}
\Omega^t\res L)=d(\widetilde{dh}\ip\mathrm{\,vol}^t)=d(*dh)=\Delta(*h),
$$
where $J$ is the complex structure and $\,\widetilde{\ }\,$ is the
isomorphism $T^*L\to TL$ set up by the induced metric on $L$. So the
equation has a solution by the Fredholm alternative and (\ref{fred}).
\end{Proof}

Thus for $\phi>0$ we consider the $L^t$s as the \emph{same} as
Lagrangian submanifolds (up to hamiltonian deformation) in the
\emph{fixed} symplectic manifold
$W^t$; it is only the SLag representative that changes as $\Omega^t$
varies. We think of this as mirror to a \emph{fixed} holomorphic bundle
in a fixed complex structure, with varying HYM connection as the mirror
K\"ahler form changes.

\begin{Lemma} \label{lem} In the analogous 2-dimensional situation of
SLags in a $K3$ or abelian surface, the obstruction does not occur.
\end{Lemma}

\begin{Proof}
Choose a real path of complex structures $W^t,\
t\in(-\epsilon,\epsilon)$ in complex structure moduli space such that
there is a nodal SLag $L^0=L_1^0\cup L_2^0$ in $W^0$. Without loss of
generality we can choose the
phase of $\Omega^t$ so that both $\omega$ and Im$\,\Omega^t$ pair to
zero on the homology class of $L^0$. Now hyperk\"ahler rotate the
complex structures so that instead the new Re$\,\Omega^t$ and
Im$\,\Omega^t$ pair to zero on the homology class of $L^0$ for all
$t$. $L^0$ is now a nodal holomorphic curve $C$ in the central $K3$.
We can understand deformations of $C$ via
deformations of the ideal sheaf $\mathcal J_C$, with obstructions in
$$
\mathrm{Ext}^2(\mathcal J_C,\mathcal J_C)\to H^{0,2}(W)\cong\C,
$$
where the arrow is the trace map and is an isomorphism by Serre
duality. Standard deformation theory shows the obstruction is
purely cohomological -- it
is the derivative of the $H^{0,2}$-component of the class
$$
[C]\in H^2(W)\cong H^{2,0}(W)\oplus H^{1,1}(W)\oplus H^{0,2}(W).
$$
But we have fixed this to remain zero by the phase condition,
so the curve deforms to all $t$ (really we should assume the family
is analytic in $t$ here and extend to $t\in\C$, or just work with
first order deformations). Hyperk\"ahler unrotating gives back a
family of SLags.
\end{Proof}

There is a notion of connect summing Lagrangian submanifolds intersecting
in a single point (probably due to Polterovich) -- see for instance
Appendix A of \cite{S1} -- which
we claim gives the smoothings $L^t$ of the singular $L^0=L_1\cup L_2$.
This follows by comparing the local models \cite{J}, \cite{S1} for
the Lagrangians; see \cite{TY} where it is studied in more detail for
a related purpose, and our conventions (slightly different from those
of \cite{S2}) are described. While
topologically we are just connect summing $L_1$ and $L_2$ by removing a
small 3-ball containing the intersection point from each and gluing the
resulting boundary $S^3$s (there are two ways, depending on orientation),
symplectically the construction does not explicitly use orientations of
the submanifolds. (Effectively we are using their relative orientation --
the canonical orientation of the sum of the tangent spaces of $L_1,\ L_2$
at the intersection point given by the symplectic structure.)

Giving $L_1$ and $L_2$ in that order produces a Lagrangian, well
defined up to hamiltonian isotopy (this will be shown in Section
\ref{Ko} in more generality; see (\ref{proj})),
$$
L_1\#L_2,
$$
with the singular union $L_1\cup L_2$ a limit point in the hamiltonian
isotopy class, which is \emph{not} itself hamiltonian isotopic
to $L_1\#L_2$ (we have seen a family of hamiltonian
deformations which has limit $L_1\cup L_2$, but the deformations
are singular at this limit).

There is also an obvious notion of \emph{graded connect sum},
which is in fact what we shall always mean by $\#$. There is a unique
grading on $L_1$ compatible with a fixed grading on $L_2$ such that
we can give a (continuous) grading to the smoothing $L_1\#L_2$. In the
case of connect summing at multiple intersection points (Section
\ref{Ko}) there is \emph{at most} one such grading; in general the
graded connect sum may not exist.

In $n$ dimensions, if $L_1$ and $L_2$ are graded such that $L_1\#L_2$
exists, then on reversing the order of the $L_i$, the graded sum that
exists is
\begin{equation} \label{+-}
L_2\#(L_1[2-n]) \quad\text{in the homology class}\quad
[L_2]+(-1)^n[L_1].
\end{equation}
Here $L[m]$ means the graded Lagrangian $L$ with its grading changed by
adding $m\pi$ to $\theta$, and the homology class of $L_1\#L_2$ is
$[L_1]+[L_2]$ using the orientations on the $L_i$s induced by the gradings.

This is closely related, as we shall see, to Joyce's obstruction,
and the lack of it in dimension 2 (Lemma \ref{lem}). In 2 dimensions,
$L_1\#L_2$ and $L_2\#L_1$ are in the same homology class, though
by a result of Seidel \cite{S1} \emph{not in general in the same hamiltonian
isotopy class},
$$
L_1\#L_2\not\approx L_2\#L_1,
$$
importantly (we use $\approx$ to denote equivalence up to hamiltonian
deformations). For $t>0$ in the above family $L^t$ is in the constant
hamiltonian deformation class of $L_1\#L_2$, for $t<0$ it is in the
different class of $L_2\#L_1$, and at $t=0$ it is $L_1\cup L_2$ -- in
neither class but in the closure of both. (For
complex $t$ the symplectic structure is no longer constant like it is
for $t\in\R$, as one can see by following through the hyperk\"ahler
rotation; thus we do not get a contradiction to the above statement
by going round $t=0$ in $\C$.) In 3 dimensions, however, the
corresponding obvious choice for a SLag on the other side of the
$\pi^t_1=0$ wall, $L_2\#L_1[-1]$, is in the wrong homology class.

In the case that the $L_i$ are Lagrangian spheres we can see this
by going round the wall
\begin{equation} \label{wall}
\phi(L_1^t)=0\simeq\phi(L_2^t),
\end{equation}
and using monodromy. In the 2-dimensional $K3$ or $T^4$ case this
works as follows.

\begin{figure}[h]
\vspace{4mm}
\center{
\input{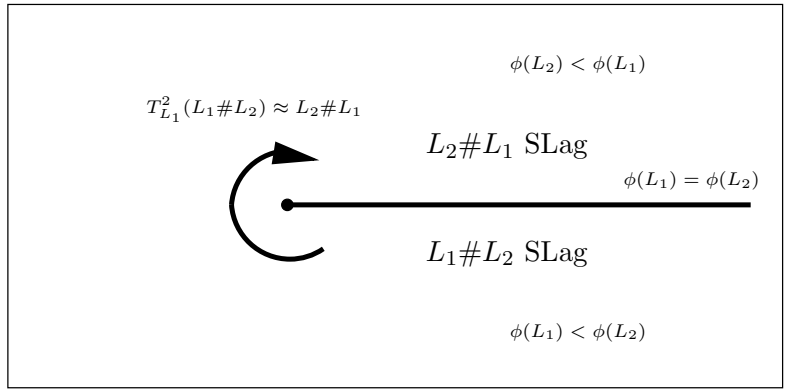}
\caption{$\left(\int_{L_1}\Omega\right)$-space, as $\Omega$
on $K3$ varies, with polar coordinates $(R,\,\phi(L_1))$}}
\vspace{2mm}
\end{figure}

Consider a disc in complex structure moduli space over which the
family of K\"ahler $K3$ surfaces (with constant K\"ahler form)
degenerates at the origin to a $K3$ with an ordinary double point
(ODP) with the Lagrangian $L_1\cong S^2$ as vanishing cycle.
A local model is the standard K\"ahler structure on $x^2+y^2+z^2=u$,
over the parameter $u$ in the unit disc in $\C$. Now base-changing
by pulling back to the double cover in $u$, $u\mapsto u^2$, we
get the 3-fold 
$$
x^2+y^2+z^2=u^2,
$$
with a 3-fold ODP which has a small
resolution at the origin putting in a holomorphic sphere resolving
the central $K3$ fibre $u=0$. Choosing a nowhere-zero holomorphic
section $\Omega_u$ of the fibrewise $(2,0)$-forms
(using the fact that the relative canonical bundle of either family is
trivial), this restricts to zero on the exceptional $\Pee^1$ (which is
homologous to the vanishing cycle $L_1$). Therefore the function
\begin{equation} \label{zero}
\int_{L_1}\Omega_u
\end{equation}
\emph{has a simple zero at} $u=0$, i.e. it vanishes to order 1 in $u$.
(The same expression vanished only as $\sqrt u$ in the original family
with the singular fibre, and as such its sign was not well defined; this
is because the class $[L_1]$ was defined globally only up to the
monodromy $T_{L_1}[L_1]=-[L_1]$, i.e. up to a sign. In our new family
the monodromy action $T^2_{L_1}$ is trivial on homology so it makes sense
to talk about the homology class $[L_1]$ in any fibre, and (\ref{zero})
is single valued.)

Then our loop of complex structures is given by taking the loop
$u=e^{it}$ and setting $\Omega^t=\Omega_{e^{it}}$. Pulling back the
K\"ahler form from the original family, we get a locally
trivial fibre bundle of symplectic manifolds over the
circle whose monodromy is the Dehn twist $T_{L_1}^2$ (since the
monodromy round the un-base-changed loop is $T_{L_1}$ \cite{S1}). As
the family no longer has a singular fibre this monodromy is trivial as
a diffeomorphism, but it is a result of \cite{S1},\,\cite{S2} that as a
symplectic automorphism \emph{it is non-trivial}. This is possible
because although the family is a locally trivial bundle of symplectic
manifolds over the punctured disc, over $u=0$ the symplectic form
becomes degenerate since it was pulled back via the resolution map.

Measuring $[L_1]$ against $\Omega_u$ as in (\ref{zero}) we see a
principle familiar in physics (in issues of `marginal stability',
and taught to me by Eric Zaslow) -- we detect a monodromy,
like the degree 1 map $S^1\to\C^{\!\times}$ given by $t\mapsto\int_{L_1}
\Omega^t$ in this example, by counting wall crossing where a certain
real part (here $\int_{L_1}\ImO^t$, or the phase $\phi_1^t$) hits
$0\simeq\phi^t_2$.

(Here we can no longer choose the phase of $\Omega$
such that $\phi^t_2=\phi(L^t_2)\equiv0$ in the whole family, as the
homology class of $L_2$ is not preserved in the family:
$$
[T_{L_1}^2L_2]=[L_2]+2[L_1].
$$
However, for a sufficiently small loop about the ODP, i.e. for
$\big|\int_{L_1}\Omega\big|$ sufficiently small, this will not
affect us much and we can write $\phi^t_2\simeq0$: we are only
interested in topological information like winding numbers and
$\phi^t_1$ crossing the wall at $\phi^t_2\simeq0$, which are
unaffected by small perturbations.)

So instead of going through the $\phi(L_1^t)=\phi(L_2^t)\simeq0$
wall we can go round it. If the loop is sufficiently
small we do not encounter any more walls where the homology class
$[L_1]+[L_2]$ can be split into classes of the same phase to possibly
make the SLag a singular union of distinct SLags of equal phase. For
instance the wall at phase 0 does not extend past $u=0$ to phase
$\phi^t_1=\pi$ (even though there $\mu^t_1=0$) --
the phase of $L_1$ is not zero but $\pi$, and is only zero for $L_1$ with
the opposite orientation, so it does not exist as a SLag (e.g. in the
hyperk\"ahler rotated situation, we are saying there is no complex curve in
$L_1$'s homology class to possibly make $L$ the nodal union of $L_1$
and something else, there is only an anti-complex curve).
So we really can go round the wall; it ends at $u=0$.

So this monodromy
description shows that on the other $t\uparrow2\pi$ side of the wall
the SLag deforming $L_2\cup L_1$ is in the hamiltonian deformation
class
\begin{equation} \label{giggs}
T_{L_1}^2L\,=\,T_{L_1}^2(L_1\#L_2)\,=\,T_{L_1}^2(T_{L_1}^{-1}L_2)\,
\approx\,T_{L_1}L_2\,\approx\,L_2\#L_1,
\end{equation}
as claimed (for the above equalities see \cite{S1}, \cite{S2}).

Notice that the alternative connect sum description of the
above Lagrangian
\begin{equation} \label{decomp}
L_2\#L_1\,=\,T_{L_1}^2(L_1\#L_2)\,\approx\,T_{L_1}^2(L_1)\#
T_{L_1}^2(L_2)\,\approx\,L_1[-2]\#T_{L_1}^2(L_2),
\end{equation}
does not violate the phase inequality to (\ref{ineq}), as
$$
-2\pi+\epsilon\approx\phi(L_1[-2])<\phi(T_{L_1}^2(L_2))\simeq0.
$$
This is why it is important here to keep track of gradings --
assigning the phase
$\epsilon$ to $\phi(T_{L_1}^2(L_1))$ would give the opposite
inequality, but one would not be able to form the above graded
connect sum without also shifting
the phase of $T_{L_1}^2(L_2)$ by $-2\pi$.

\begin{figure}[h]
\vspace{4mm}
\center{
\input{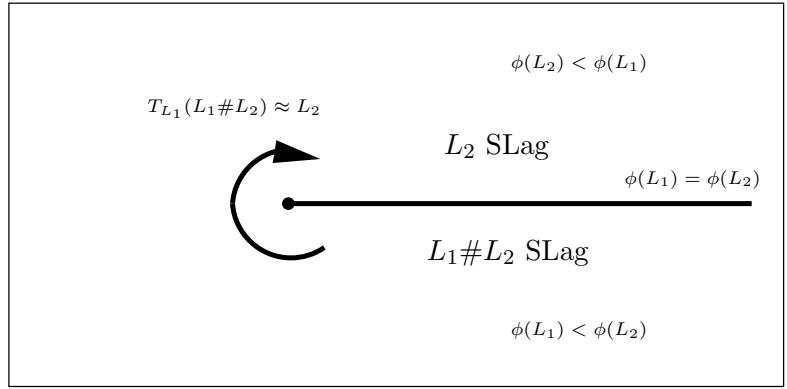}
\caption{$\left(\int_{L_1}\Omega\right)$-space, as $\Omega$
on a 3-fold varies, with polar coordinates $(R,\,\phi(L_1))$}}
\end{figure}

The 3-fold case (which Dominic Joyce has also, independently, considered)
is slightly different; we need only take a single Dehn twist
$T_{L_1}$ corresponding to the local family
$$
x^2+y^2+z^2=u,
$$
over $u\in\C$ to get a winding number one loop in the phase of $L_1$.
This is because
$$
T_{L_1}L_1\approx L_1[1-n]
$$
in dimension $n$, so in 3 dimensions the homology class $[L_1]$ is
preserved instead of being reversed. The corresponding picture is
displayed in Figure 2.

Again there is a SLag on the other side of the $\phi=0$ wall, but
it is in the wrong homology class $[L_2]$:
\begin{equation} \label{cole}
T_{L_1}L\approx L_2.
\end{equation}
Analogously to (\ref{decomp}) this has a number of decompositions
as connect sums induced by monodromy,
$$
T_{L_1}(L_1\#L_2)\,\approx\,L_1[-2]\#(L_2\#(L_1[-1]))\,\approx\,L_2
\,\approx\,(L_1\#L_2)\#(L_1[\,1\,]),
$$
none of which violate the phase inequality (\ref{ineq}). The only other
obvious choice for a (S)Lag on the other side of the $\phi=0$ wall
(given the $K3$ result) is $T_{L_1}^2(L_1\#L_2)\approx L_2\#(L_1[-1])$;
this however is also in the wrong homology class, and in any case
does violate (\ref{ineq}) and so, by Joyce's analysis, should not
be represented by a SLag. Thinking of $T_{L_1}^2$ as rotating through
$-4\pi$ in Figure 2, it is at roughly $-3\pi$ that the phase inequality
(\ref{ineq}) gets violated, and the $-\pi$ rotation of $L_2$ splits
as a SLag into the union of the $-\pi$ rotations of $(L_1\#L_2)$ and
$L_1[\,1\,]$: these both have phase approximately zero.

\subsection*{A holomorphic bundle example}

These phenomena are similar to wall-crossing in bundle theory on the
complex side -- in a real one-parameter family of K\"ahler forms, for
fixed complex structure, stable holomorphic bundles for $t>0$ can become
semistable at $t=0$ and unstable for $t<0$.

An example that mirrors Joyce's is the following. Suppose we
have two stable bundles (or coherent sheaves) $E_1$ and $E_2$ with
$$
\mathrm{Ext}^1(E_2,E_1)\cong\C.
$$
This is $H^1(E_1\otimes E_2^*)$ in the case of
bundles and is the mirror \cite{K} of the one dimensional
Floer cohomology $HF^*(L_2, L_1)\cong\C$ that is defined by the single
intersection point of $L_1$ and $L_2$ (see Section \ref{Ko} for more
details of this, and an explanation of why we are dealing with Ext$^1$
and $HF^1$ here). We then form $E$ from this extension class
\begin{equation} \label{E}
0\to E_1\to E\to E_2\to0.
\end{equation}
Take a family of K\"ahler forms $\omega^t$ such that
$\mu^t(E_2)-\mu^t(E_1)$ is the same sign as $t$ (here $\mu^t(F)=
c_1(F)\,.\,(\omega^t)^{n-1}/$\,rk\,$(F)$ is the slope of $F$ with
respect to $\omega^t$). Supposing that the $E_i$ are stable for all
$t\in(-\epsilon,\epsilon)$, we claim that $E$ is stable
for sufficiently small $t>0$, while it is destabilised by $E_1$ for
$t\le0$.  Without loss of generality take $\mu^t(E_2)=\mu$ fixed,
and $\mu^t(E_1)=\mu-t$. As $E_2$ is stable, for $t$ sufficiently small
there are no subsheaves of $E_2$ of slope greater than $\mu-t$, so for
any stable destabilising subsheaf $F$ of $E$, the composition
$$
F\into E\to E_2
$$
cannot be an injection (unless it is an isomorphism, but (\ref{E}) does
not split. So $F\cap E_1\ne0$, and the quotient $Q=F/(F\cap E_1)$ has
slope $\mu(Q)>\mu(F)>\mu-t$ by the stability of $F$ and
instability of $E$. But $Q$ injects into $E_2$, which we know is
impossible.

In the 2-dimensional case, by Serre duality Ext$^1(E_1,E_2)
\cong\,$Ext$^1(E_2,E_1)^*\cong\C$on $K3$ or $T^4$, so for $t<0$ we
can instead form an extension
\begin{equation} \label{ext}
0\to E_2\to E'\to E_1\to0,
\end{equation}
to give a new bundle $E'$ which is also stable, and has the same Mukai
vector
$$
v(E')=v(E_1)+v(E_2);
$$
compare (\ref{+-}). At $t=0$ we take the (polystable)
bundle
$$
E_1\oplus E_2.
$$
This is because the semistable extension (\ref{E}) no longer admits a
Hermitian-Yang-Mills metric, but $E_1\oplus E_2$ does. Also, the algebraic
geometry of the moduli problem shows that while a semistable bundle gets
identified in the moduli space with the other (``S-equivalent'') sheaves
in the closure of its gauge group orbit, there is a distinguished
representative of its
equivalence class -- the polystable direct sum (of the Jordan-H\"older
filtration, which here is $E_1\oplus E_2$).

Thus, while the HYM connections vary, the bundle
has only 3 different holomorphic structures -- for $t>0,\ t=0,$ and $t<0$.
Put another way (to spell out the analogy with the Lagrangians $L^t,\ 
L_1,\ L_2$) as $\omega_t$ varies with $t>0$ we take different points in
a fixed complexified gauge group orbit, and at $t=0$ we take as limit point
something in a different orbit that is nonetheless in the closure of the
$t>0$ (and $t<0$) orbit. The \emph{stable} deformations of the polystable
$E_1\oplus E_2$ (which we are thinking of as the mirror of the singular
union $L_1\cup L_2$, of course) are precisely (\ref{E}) for $t>0$ and
(\ref{ext}) for $t<0$.

In the 3-fold case, however, Serre duality gives Ext$^2(E_1,E_2)\cong\,
$Ext$^1(E_2,E_1)^*\cong\C$ instead, and so no stable extension (\ref{ext}).
In fact one would expect there to be no stable bundle with the right Chern
classes; instead the one dimensional Ext$^2$ gives us a complex $E'$ in
the derived category $D^b(M)$ fitting into an exact sequence of complexes
$$
0\to E_2\to E'\to E_1[-1]\to0,
$$
where $E_1[-1]$ is $E_1$ shifted in degree by one place to the right as
a complex. This has Mukai vector
$$
v(E')=v(E_2)-v(E_1),
$$
compare (\ref{+-}). Thus, just as in the
case of SLags, as we pass through $t=0$ there is no natural stable
object on the other side in the same homology class in 3 dimensions
(though there is in 2 dimensions) and so an element of the appropriate
moduli space disappears.

In fact, as in the Lagrangian example, the natural stable object
on the other side of the wall is $E_2$ if we consider monodromy.
The mirror of the symplectic Dehn twists of above are described
in \cite{ST} (in the case that the bundles $E_i$ are \emph{spherical}
in the sense of \cite{ST}: Ext$^k(E_i,E_i)\cong H^k(S^n;\C)$; this is
the natural mirror analogue of the
$L_i$s being spheres). These are the twists $T_{E_1}$ of \cite{ST}
on the derived category of the Calabi-Yau that act on the extension
bundle $E$ of (\ref{E}) to give precisely the extension (\ref{ext}),
$$
T_{E_1}^2E=E'
$$
(compare (\ref{giggs})), as a short calculation using \cite{ST}
shows. Similarly
$$
T_{E_1}E=E_2,
$$
the analogue of (\ref{cole}). (In both of these calculations it is
important to calculate this monodromy in the derived
category; in the $K3$ case the action of $T_{E_1}^2$ is trivial on
K-theory and cohomology, and we cannot distinguish between (\ref{E})
and (\ref{ext}), but they are very different as holomorphic bundles
and as elements of the derived category.)

The mirror wall crossing, with a SLag splitting into two and then
disappearing, is interpreted in \cite{DFR} (and in \cite{SV} in a
different case) as the
state it represents decaying as we reach a point of `marginal
stability'. Despite this dealing with only SLags (and so with only
a priori \emph{stable} Lagrangians in our mathematical sense of
stability), this suggestive language does in fact have something to say
about the stability, in our sense of group actions, of (non-special)
Lagrangians, by considering the nodal limit $L_1\cup L_2$ to be a
\emph{semistable} Lagrangian.

Thus the Lagrangian $L_1\#L_2$ (which
\emph{always exists as a Lagrangian} as the complex
structure varies with fixed K\"ahler form) becomes semistable
at $t=0$ and is represented by something in a different orbit of the
hamiltonian deformation symmetry group (but in the closure of the
original orbit), and is unstable for $t<0$ so exists there only as a
Lagrangian and \emph{not} as a SLag. This, and the bundle analogue
described above, leads us to think of the Lagrangian $L_1$
as destabilising $L=L_1\#L_2$ when $\phi(L_1)\ge\phi(L_2)$. This
motivates the now obvious definition of stability in Section \ref{st};
first we explain more about the connections to mirror symmetry, and
generalisations to connect sums at more intersection points.

\section{Relationship to Kontsevich's mirror conjecture} \label{Ko}

The inspiration behind most of this paper is of course Kontsevich's
mirror conjecture \cite{K}. In particular, Kontsevich proposes that
the graded vector spaces Ext$^*$ and $HF^*$ should be isomorphic for
mirror choices of bundles $E_i$ and graded Lagrangians $L_i$
(or more exotic objects in their derived categories)
$$
HF^*(L_2,L_1)\cong\,\mathrm{Ext}^*(E_2,E_1);
$$
this corresponds to the equality of (graded) morphisms on both sides.
Here $HF^*$ is Floer cohomology \cite{Fl} -- a symplectic refinement
of the intersection number of $L_1$ and $L_2$ -- which can be $\mathbb
Z$-graded for graded Lagrangians \cite{S2}, whenever it is defined
\cite{FO3}, \cite{Fu1}. (More precisely it is the cohomology of a chain
complex built out of the free vector space generated by the intersection
points, with the differential defined by counting holomorphic
discs with boundary in the Lagrangians running from one intersection
point to another.) In mirror symmetry, and so in this paper, one
should only really consider those Lagrangians whose Floer cohomology is
well defined \cite{Fu1}.

Thus the point of intersection of the $L_1$ and $L_2$ of the last
section define the Floer cohomology $HF^*(L_2,L_1)\cong\C$, and the
grading of \cite{S2} is designed specifically so that $L_1\#L_2$ can
be graded precisely when the relative gradings of the $L_i$ force
the Floer cohomology to be concentrated in degree 1; $HF^*(L_2,L_1)=
HF^1(L_2,L_1)$. We then think of the connect sum $L_1\#L_2$ as being
mirror to the extension (\ref{E}) defined by Ext$^1(E_2,E_1)\cong\C$.
Fukaya, Seidel, and perhaps others have also proposed that Lagrangian
connect sum should be mirror to extensions \cite{Fu2}, \cite{S3}.

We also consider connect sums of Lagrangians intersecting at $n$
points $p_i$. Then the connect sum is not unique up to hamiltonian
deformation: $H^1$ is added to the Lagrangian as loops between the
intersection points, giving additional deformations of its hamiltonian
isotopy class. The upshot is that there is a scaling of the neck of
the connect sum at each intersection point; we denote any such
resulting Lagrangian by $L_1\#L_2$. Since we insist on all intersection
points having Floer (Maslov) index one (so that the connect sum can be
graded), the Floer differential vanishes in this case, and
these scalings define a class in $HF^1(L_2,L_1)$.

Deformations (up to those which are hamiltonian) as such a connect sum
are given by the elements of
$$
H^1(L_1\#L_2)\cong H_{n-1}(L_1\#L_2)
$$
spanned by the $S^{n-1}$ vanishing cycles $S_i$ at the points of
intersection $p_i\in L_1\cap L_2$. Given a particular connect sum,
the deformation represented by $\sum_ia_iS_i$ simply scales the local
gluing parameter in a Darboux chart around each $p_i$ by a factor
$(1+a_i)$ (here $a_i$ is considered to be infinitesimal). Since
the sum of these spheres
separates $L_1\#L_2$ into $L_1\backslash\cup\{p_i\}$ and
$L_2\backslash\cup\{p_i\}$ and so is zero in homology
$$
\sum_i[S_i]=\pm\partial[L_1\backslash\cup\{p_i\}]=
\mp\partial[L_2\backslash\cup\{p_i\}]=0\in H_{n-1}(L_1\#L_2),
$$
the infinitesimal deformation represented by $\sum_iS_i$ is zero
(it is pure hamiltonian) and dividing out gives the projectivisation
\begin{equation} \label{proj}
\Pee(\oplus_i\R_{p_i}).
\end{equation}
(Replace $\R$ by $\C$ when including flat bundles and their gluing
parameters at the $p_i$s.)  This
explains the earlier claim that connect sums at one intersection
point are uniquely defined up to hamiltonian deformations.
More precisely, when holomorphic discs are
taken into account and we consider only those Lagrangians whose Floer
cohomology is defined \cite{FO3}, hamiltonian deformation classes of
connect sums \emph{whose Floer cohomology can be defined}
should be parameterised by $\Pee(HF^1(L_2,L_1))$.
(On the mirror side \emph{isomorphism classes} of extensions of $E_2$
by $E_1$ are parametrised by $\Pee\mathrm{\,Ext}^1(E_2,E_1)$.)

We would then expect that the resulting connect sum has a canonical
homomorphism from $L_1$; that is there should be a canonical element
$$
\id_{L_1}\in HF^0(L_1,L_1\#L_2)
$$
for any graded Lagrangians $L_i$ for which the graded connect sum
exists. While a local model suggests this is true (see for instance
\cite{TY}), a complete proof is
still not available. This homomorphism we think of as expressing $L_1$
as a \emph{subobject} of $L_1\#L_2$; i.e. as giving an injection.
It should be emphasised that subobject does not make sense in a
triangulated category such as the derived Fukaya category of Lagrangians;
in the context of the derived category of sheaves, subobject
only makes sense for an abelian category such as that of the sheaves
themselves (i.e. complexes with cohomology in degree zero only).
What we are proposing is that it also makes sense in the category of
(complexes of sheaves mirror to) graded Lagrangians, and is vital to
make definitions of stability (which involve such subobjects). While
there are now more Homs to consider, in particular those of higher
order (i.e. Homs to Lagrangians shifted in phase by some $2\pi n$),
the targets of these Homs have higher phase and so do not disturb the
definition of stability below -- this is seemingly a huge piece of luck
that means we can extend the stability condition for bundles to all
Lagrangians. For similar reasons, the many
connect sum decompositions of the $L_i$s given in the last section
also do not destabilise them.

There are other operations, however, which can also be thought of
as Ext$^1$-type extensions. For instance, taking the product of a
single Lagrangian curve $L_1$ in $T^2$ with a (graded) connect sum
$L_2\#L_3$ in another $T^2$, we get a
Lagrangian $L_1\times(L_2\#L_3)$ in $T^4$ which is some kind of
extension of the Lagrangians $L_1\times L_2$ and $L_1\times L_3$ in
$T^4$. Supposing that the $L_i$s are mirror to some (complexes of)
sheaves $E_i$, and that the connect sum $L_2\#L_3$ is mirror to an
extension represented by an element $e\in$\,Ext$^1(E_3,E_2)$. Then
by the K\"unneth formula for sheaf cohomology, we see that
$L_1\times(L_2\#L_3)$ is indeed mirror to an extension
$$
\id\,\otimes e\in\mathrm{Hom\,}(E_1,E_1)\otimes\mathrm{Ext}^1
(E_3,E_2)=\mathrm{Ext}^1(E_1\boxtimes E_3,E_1\boxtimes E_2),
$$
and so this sort of relative connect sum (which is not $\#$ on
$T^4$: $L_1\times L_2$ and $L_1\times L_3$ do not intersect
transversely) should also be considered.

So we consider Lagrangians $L_1,\,L_2$ intersecting \emph{cleanly}
(see e.g. \cite{S1} Definition 2.1), that is $N=L_1\cap L_2$ is a smooth
submanifold, and $TN=TL_1\res N\cap TL_2\res N$. Basic results of
Weinstein allow us to identify a neighbourhood of $N$ with a
neighbourhood of the zero section $N$ in $T^*N\oplus E$, where the total
space of $T^*N$ has its canonical symplectic structure, and
$$
E\equiv(TL_1\res N)/TN\,\oplus\,(TL_2\res N)/TN
$$
is the annihilator, under the symplectic form, of $TN\subset TX\res N$
(to which the symplectic form therefore descends, making $E$ a
symplectic bundle).

Choosing a metric on $E$, compatible with its symplectic structure,
such that its transverse
subbundles $(TL_1\res N)/TN,\ (TL_2\res N)/TN$ are orthogonal, we can now
perform the family connect sum of these, over the base $N$, since
the local model in \cite{S1} is $O(n)$ invariant. As before we insist
that this can be compatibly graded again denote it by $\#$; given a
grading on $L_1$ there will be at most one grading on $L_2$ such that
this graded relative connect sum exists.

It should be noted that although such a clean intersection could be
hamiltonian isotoped to be transverse, the resulting intersection points
would not necessarily all be of Floer/Maslov index one, and so the
pointwise graded connect sum could not be formed at every point; we
would end up with an immersed Lagrangian. Studying which immersed
Lagrangians should be included in the Fukaya category, and which
embedded Lagrangians they should be considered equivalent to, is
an important part of mirror symmetry and will need to be better understood
to refine our conjecture. For instance forming extensions of bundles
which also have nonzero homorphisms between them would appear to
be mirror to forming connect sums between graded Lagrangians at index
one intersection points, leaving the index zero intersection points
immersed. In general one would like to consider two objects of the Fukaya
category to be equivalent if their Floer cohomologies with any other
objects are the same. This would include hamiltonian deformation
equivalence, but also more exotic equivalences for immersed
Lagrangians (thanks to Paul Seidel for pointing this out to me).
A start in understanding the Floer cohomology of immersed
Lagrangians is \cite{Ak}; in the present paper we are largely ignoring
singularities.

\section{Stability} \label{st}

\begin{Definition}
Take graded Lagrangians $(L_1,\theta_1)$ and $(L_2,
\theta_2)$, hamiltonian isotoped to intersect cleanly, and such
that the graded (relative) Lagrangian connect sums
$(L_1\#L_2,\theta_1\#\theta_2)$ exist as above. Then a
Lagrangian $L$ of Maslov class zero is said to be
destabilised by the $L_i$ if it is hamiltonian isotopic to such an
$L_1\#L_2$, and the phases (\emph{real} numbers, induced by the
gradings) satisfy
$$
\phi(L_1)\ge\phi(L_2).
$$
If $L$ is not destabilised by any such $L_i$ then it is called stable.
\end{Definition}

\begin{Remarks}
\begin{description} \item[\ \ $\bullet$\ \ ]
There is an obvious notion of a flux homomorphism for isotopies
of smooth Lagrangians, taking a deformation to an element of
$H^1(L;\R)$ (and linearising to give the usual deformation theory of
Lagrangians). Namely, take a deformation $\Phi_t(L)$ through a
vector field $X_t,\ t\in[0,1]$ to the one form
$$
\int_0^1(X_t\ip\omega)dt\in H^1(L;\R).
$$
Alternatively, the homomorphism takes a loop $\gamma\subset L$,
tracing out the 2-cycle $f(\gamma\times[0,1])$ in $W$ under the
isotopy, to the real number $\int_{\gamma\times[0,1]}\omega$.
(See Chapter 10 of \cite{MS} for the analogous map for
symplectomorphisms.) If the isotopy $\Phi_t$ is hamiltonian, the
flux is zero; the converse is also easily proved using the methods of
(\cite{MS} Theorem 10.12): we may assume without loss of generality
that the 1-form $\int_0^1X_t\ip\omega$ is identically zero in
$\Omega^1(L)$. [To see this, write the 1-form as
$d\phi$, and compose the deformation with
the time one map of the hamiltonian flow with vector field
$d\phi\ip\omega^{-1}$; this does not alter the flux in $H^1(L;\R)$ or
the property of being hamiltonian.] Then let $\Sigma^s$ be the closed
1-form on $L$ defined by
$$
\Sigma^s=\int_0^sX_t\ip\omega dt,
$$
and let $\Psi^s_t$ be the corresponding flow through time $t$. Then
the flow $\phi_t=\Psi^t_t\circ\Phi_t$ is the corresponding
\emph{hamiltonian} flow from $\Phi_0(L)$ to $\Phi_1(L)$; see
\cite{MS}.

Thus it is not too hard to check if two Lagrangians
are hamiltonian deformations of each other, at least through
smooth Lagrangians, if we know they are deformations of each
other as Lagrangians. This second condition, however, is harder
to test, as the results of \cite{S1} demonstrate.

\item[\ \ $\bullet$\ \ ] As mentioned in the last section,
holomorphic discs are
crucial in both mirror symmetry and Floer cohomology; thus one
should perhaps restrict attention in the above definition to
those Lagrangians whose Floer cohomology is defined \cite{FO3},

\item[\ \ $\bullet$\ \ ]As pointed out to me by Conan Leung, this
definition and the resulting conjecture below may only be reasonable
close to the
large complex structure limit point where the mirror symmetric
arguments used to motivate the conjecture are most valid.
\end{description}
\end{Remarks}

\begin{Conj}
A Lagrangian of Maslov class zero has a special Lagrangian in its
hamiltonian deformation class if and only if it is stable, and this
SLag representative is unique.
\end{Conj}

Again, we have been vague about singularities: which we allow, and what
hamiltonian deformation equivalence would mean for them. We might also
want to restrict attention to those Lagrangians whose Floer cohomology
exists \cite{FO3}, and whose Oh spectral sequence $H^*(L)\Rightarrow
HF^*(L,L)$ \cite{Oh} degenerates; this will be discussed more in \cite{TY}.
We might also want to restrict to Lagrangians whose phase function
varies only by a certain bounded amount; in the example worked out in
\cite{TY} this is required. In \cite{TY}
it is shown there that the gradient of the
norm-squared $|m|^2$ of the moment map can be taken to be the mean
curvature vector of the Lagrangian, so mean curvature flow (which
is hamiltonian for Maslov class zero) should
converge to this SLag representative if the Lagrangian is stable
and the phase satisfies certain bounds.

\section{The 2-torus}

Everything works rather simply on $T^2$; Grayson \cite{G}, building
on work of Gage, Hamilton and others (e.g. \cite{GH}), has shown that
mean curvature flow for curves (of Maslov class zero) converges to
straight lines and so we get the mirror symmetric analogue 
of Atiyah's classification \cite{At} of sheaves on an elliptic curve
-- they are basically all sums of stable sheaves. The only exceptions
are the non-trivial extensions of certain
sheaves by themselves; these correspond to thickenings of the
corresponding special Lagrangian (giving fat SLags, as they are known
in Britain, or multiply-wrapped cycles in physics speak).

\begin{figure}[ht]
\vspace{1mm}
\center{\input{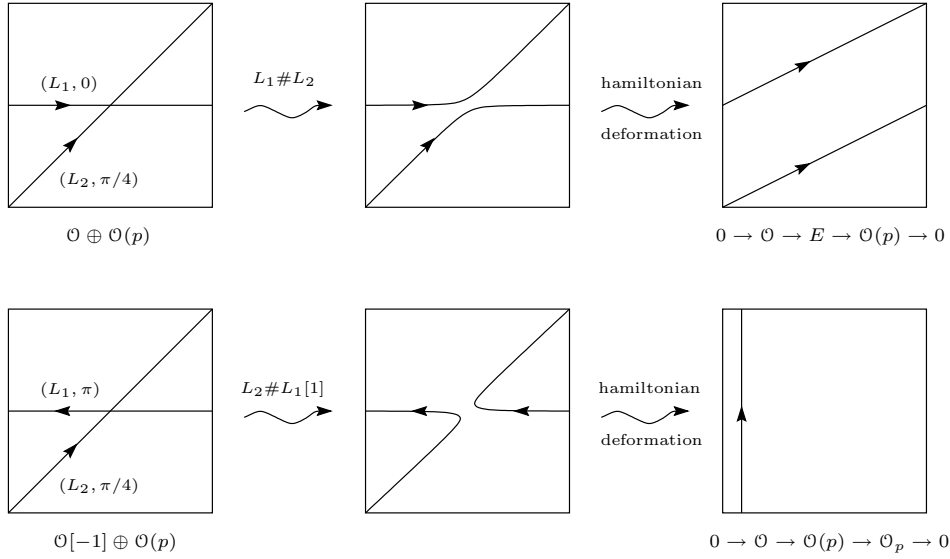}}
\caption{$L_1\# L_2$ and $L_2\#(L_1[1])$, equivalent SLags, and
their mirror sheaves}
\vspace{1mm}
\end{figure}

We give an example to demonstrate why one cannot form smooth
\emph{unstable} Lagrangians on $T^2$ in Figure 3.
First, giving $L_1$ and $L_2$ the gradings such that their
phases are 0 and $\pi/4$, we expect $L_1\#L_2$ to be stable, and
indeed we see it is hamiltonian deformation equivalent to the slope
$1/2$ SLag mirror to the stable extension $E$ of $\OO$ by $\OO(p)$
(where $p$ is a basepoint of $T^2$ with corresponding line bundle mirror
to the diagonal SLag drawn).

If one then tries to form an unstable SLag $L_2\#L_1$, the \emph{graded}
connect sum does not exist -- the phase would become
discontinuous. To form $L_2\#L_1$ we see from the diagram that we have to
take the phase of $L_1$ to be $\pi$, thus reversing its orientation, and
in fact forming $L_2\#(L_1[1])$.
Then the stability inequality (\ref{ineq}) is not violated, and in fact
this Lagrangian is stable and hamiltonian deformation equivalent to
the SLag in $T^2$ represented by the vertical edge of the square (and so
drawn with a little artistic license in Figure 3). Under the mirror
map this corresponds to replacing the extension Ext$^1$ class by a Hom
(as we have shifted complexes of sheaves by one place) and taking the
cone of this in the derived
category; this is the cokernel $\OO_p$ of Figure 3.

As pointed out to me by Markarian and Polishchuk, one can play with
lots of pictures of connect sums on tori to recover descriptions
of certain moduli of sheaves, their special cycles (for instance
where one connect-sum neck parameter goes to zero), and so forth,
giving results similar to some of those in \cite{FO}.

This example can be extended to show that we cannot form the graded
connect sum $L_1\#L_2$ of any two Lagrangians (via a class in
$HF^*(L_1,L_2)$) if $\phi(L_1)>\phi(L_2)$. Namely, replace
$L_1$ and $L_2$ by their hamiltonian deformation equivalent SLag
representatives, which are straight lines of constant phase $\theta_i=
\phi(L_i)$. As Figure 3 shows, $L_1\#L_2$ can be compatibly
graded about an intersection point if and only if we have the local
inequalities
$$
\theta_2>\theta_1>\theta_2-\pi.
$$
Thus we require $\phi(L_2)>\phi(L_1)$. (We will explain this kind
of phenomenon more generally in \cite{TY} in terms of the grading
on Floer cohomology.) Each intersection point is
Floer coclosed since the Floer grading is the same
as the relative orientation of the Lagrangians, mod 2, and this
is the same at each intersection point of the straight lines. So each
possible connect sum of the SLags defines a class in $HF^*$, and
any other connect sum, defined on hamiltonian deformations of $L_1$
and $L_2$ by a class in $HF^*$, will be hamiltonian deformation
equivalent to the appropriate connect sum of the SLags, and so
satisfy the same phase inequality.

If two smooth Lagrangians have the same phase then their
representative SLags will either be the same or disjoint parallel
SLags. Either way there are no connect sums (though as mentioned above
to account for the mirror symmetry of bundles one should also include
non-trivial thickenings of SLags in the Fukaya category).

So unstable Lagrangians do not exist, and by the result of \cite{G}
mentioned earlier, the conjecture is true on $T^2$.

Thus complex dimension 1 is too simple -- in trying to make the phase
of one Lagrangian become larger than the phase of another, the two must
cross, thus reversing their relative orientations and changing the
order of the connect sum. Far more complicated phenomena arise
in 2 and 3 dimensions, however.

\noindent {\tt thomas@math.harvard.edu} \newline
\noindent Department of Mathematics, Harvard University, One
Oxford Street, Cambridge MA 02138. USA.

\end{document}